\documentclass[11pt]{amsart}
\usepackage{amssymb,mathtools}
\usepackage[hidelinks]{hyperref}
\usepackage{xcolor}
\newtheorem{theorem}{Theorem}[section]
\newtheorem{corollary}[theorem]{Corollary}
\newtheorem{proposition}[theorem]{Proposition}
\newtheorem{lemma}[theorem]{Lemma}
\theoremstyle{definition}
\newtheorem{definition}[theorem]{Definition}
\theoremstyle{remark}
\newtheorem*{classicaltheorem}{Theorem}
\newtheorem{remark}[theorem]{Remark}
\newtheorem{example}[theorem]{Example}
\newcommand{\U}{\mathcal U}
\newcommand{\PP}{\mathbb P}
\newcommand{\OO}{\mathcal O}
\newcommand{\mm}{\mathfrak m}
\newcommand{\pp}{\mathfrak p}
\newcommand{\qq}{\mathfrak q}
\newcommand{\Frac}{\operatorname{Frac}}
\newcommand{\res}{\operatorname{res}}
\newcommand{\NN}{\mathbb N}
\newcommand{\RR}{\mathbb R}
\definecolor{revisionpurple}{RGB}{108,43,145}
\newif\ifshowrevisions
\showrevisionstrue

\title{The projective fiber of an indeterminate differential quotient}
\author{Wei Li}
\address{State Key Laboratory of Mathematical Sciences, Academy of Mathematics and Systems Science, Chinese Academy of Sciences, Beijing, China}
\email{liwei@mmrc.iss.ac.cn}
\subjclass[2020]{12H05}
\keywords{Differential specialization, indeterminate differential quotient,
monotone valuation, projective fiber}
\thanks{ The author was partially supported by the Strategic Priority Research Program of the Chinese Academy of Sciences (grant XDA0480503), the National Key R\&D Program of China (grant 2023YFA1009401), and the National Natural Science Foundation of China (grants 12671568 and 12288201). }
\date{}
\begin{document}
\begin{abstract}
 We prove Ritt's conjecture on the values admitted by a differential
quotient at a common zero of its numerator and denominator. The
projective fiber is either a singleton or the entire projective line.
The result holds in characteristic zero for finitely many unknowns
and commuting derivations. In the singleton case, the lowest nonzero
homogeneous pair in the translated numerator and denominator represents
the unique value. The proof constructs, for each nonzero prime
differential ideal of a differential affine domain, a monotone real
valuation on the function field of a strictly smaller prime with
center at the prescribed prime. A comparison of weighted perturbations
along this valuation gives uniqueness whenever a graph equation
specializes to a nonzero differential polynomial in the ordinate. 
\end{abstract}

\maketitle

\section{Introduction}

 Let $K$ be a differential field of characteristic zero with a finite
family $\Delta=\{\delta_1,\ldots,\delta_m\}$ of commuting derivations,
and let $\U$ be a universal differential extension of $K$. Let $F,G$
be nonzero, relatively prime differential polynomials in
$Y=(Y_1,\ldots,Y_n)$ over $K$. Which values does the indeterminate quotient $F/G$
admit at a common zero $\bar{a}$ of $F$ and $G$? 
In his 1932 book \emph{Differential Equations from the Algebraic Standpoint},
Ritt conjectured that in the ordinary differential case a \emph{form quotient}
admits either a unique value or every value at a common zero of its
numerator and denominator \cite[p.~61]{ritt1932}.
He later proved the conjecture for one derivation and one unknown
when $FG$ has order one \cite{ritt}. We prove the dichotomy at every
point of $\U^n$. 

 Kolchin stated the conjecture
for finitely many unknowns and commuting derivations  in his ICM address,  taking the origin
as the base point. He presented this formulation as a special case
of the problem of extending differential specializations
\cite[p.~272]{icm}. For one derivation, Blum characterized
the admitted values, which he called \emph{adherence values}, by
specialization of the generic graph point. He recorded Kolchin's
proof that the fiber at the origin is nonempty for any number of
unknowns, and also observed that every projective value is admitted
if the numerator and denominator have the same lowest total degree
and the ratio of their homogeneous parts of that degree lies outside
the base field  \cite{blum}. For one unknown,
Blum gave sufficient conditions for uniqueness and for fullness in
terms of the components of the graph equation. In particular, he
recovered Ritt's first-order result by an argument due to Kolchin
\cite[Theorem~1]{blum}. Remark~\ref{rem:blum-order} answers his
question about the order hypothesis in part~(iii) of that theorem.

  Choose $t\in\U^n$ differentially algebraically independent over $K$
and put $u=F(t)/G(t)$. The \emph{projective fiber} over
$b\in\U^n$ consists of the finite values $z\in\U$ for which
$(t,u)$ specializes over $K$ to $(b,z)$, together with
infinity if $(t,1/u)$ specializes over $K$ to $(b,0)$.
Definition~\ref{def:projective-fiber} describes this fiber in terms
of the Kolchin closure of the graph.  Now we state our main result.

\begin{theorem}\label{thm:main}
 For every $b\in\U^n$, the projective fiber of $F/G$ over $b$
is either a singleton or $\PP^1(\U)$. To describe the singleton,
let $E=K\langle b\rangle$ and take an $n$-tuple of differential indeterminates
$T$ over $E$.   Write $F_d, G_d$ for the homogeneous parts of $F(b+T),G(b+T)$ in
the least degree $d$ in which at least one of them is nonzero.
If the fiber is a singleton, there exists a pair
$(p,q)\in E^2\setminus\{(0,0)\}$ such that
\(                     qF_d=pG_d\,\text{in } E\{T\},
\)
and the fiber consists of the unique point $[p:q]\in\PP^1(E)$.
\end{theorem}

We record the case $b=0$ of Theorem~\ref{thm:main},
which gives Ritt's conjecture in the form stated by Kolchin in   \cite[\S2, p.~272]{icm}.
\begin{corollary}\label{cor:ritt}
 Suppose $F(0)=G(0)=0$. If $F/G$ admits two distinct values at
the origin, it admits every point of $\PP^1(\U)$ there. Otherwise it
admits exactly one value, and this value belongs to $\PP^1(K)$. 
\end{corollary}

  Two examples illustrate the difficulty when the lowest
homogeneous parts do represent a point over the ground field.
With one derivation, both $X^2/Y$ and $(Y')^2/Y$ have leading
projective direction $[0:1]$ at the origin. The first quotient
admits every value (Remark~\ref{rem:change-degree}), whereas the
second admits only zero \cite[p.~681]{blum}. Thus even when
Kolchin's constant-parameter substitution gives a value in $K$,
one must determine what other specializations can yield.

Our proof combines Kolchin's perturbation with descent on
the prime differential ideals of the graph. We construct monotone
valuations with prescribed centers to specialize from a strictly
smaller prime to a given one. Suppose the perturbation value at $a$
is defined over the valued field. If it failed to reduce to the value
at the residue point $b$, some weight for the parameter would give
the quotient a differentially transcendental residue. A graph equation
that remains nonzero in the ordinate at $b$ rules this out, and descent
then forces every admitted value to agree with the one given by the perturbation. If all graph
equations specialize identically to zero at $b$, the fiber is the
whole projective line.

The paper is organized as follows.  Section~\ref{sec:prelims} introduces differential specialization,
projective fibers, and the requisite valuation theory.
Section~\ref{sec:valuation} constructs monotone valuations with
prescribed centers, and Section~\ref{sec:comparison} uses them to
compare radial values under specialization. Section~\ref{sec:uniqueness}
passes from uniqueness in one affine chart to Theorem~\ref{thm:main}
and its corollary. 

\vskip5pt
 { \bf AI disclosure}.  
 The author acknowledges the assistance of  GPT-6 Astra during the research leading to this paper.
A proof of Theorem 1.1 in its differential-specialization formulation has been formalized in Lean 4 by Codex.
The author developed and verified the final mathematical arguments, and takes responsibility for the mathematical content and final form.  

\section{Preliminaries}\label{sec:prelims}

 In this section, we give some basic notions and preliminary results for differential specialization, projective fibers
and valuations.  For background on differential algebra we refer to Ritt and Kolchin
\cite{ritt1950,kolchin} and for valuation theory to Engler and Prestel
\cite{engler-prestel}. 

\subsection[Differential specialization]{ {Differential specialization}}
\label{subsec:specialization}

Let $(K,\Delta)$ be a differential field of characteristic zero,
where $\Delta=\{\delta_1,\ldots,\delta_m\}$ is a finite family
of commuting derivations, and fix a universal differential
extension $\U$ of $K$. 
An element $c\in\U$ is a \emph{constant} if $\delta_i c=0$ for
every $i$. We take $0\in\NN$ and, for $\nu\in\NN^m$, write
\[
    \delta^\nu=\delta_1^{\nu_1}\cdots\delta_m^{\nu_m},
    \qquad |\nu|=\sum_i\nu_i,\qquad \nu!=\prod_i\nu_i!.
\]
The differential polynomial ring
\[
    K\{Y\}=K[\delta^\nu Y_j:\nu\in\NN^m,\ 1\le j\le n]
\]
is the ordinary polynomial ring in the jets $\delta^\nu Y_j$. 
The
derivation $\delta_i$ sends $\delta^\nu Y_j$ to
$\delta^{\nu+e_i}Y_j$, where $e_i$ is the $i$th standard basis
vector. A monomial in  the $\delta^\nu Y_j$ is called a differential monomial in $Y$. 
For a tuple $a$ in a differential field extension of $K$,
write $K\{a\}$ and $K\langle a\rangle$ for the differential
$K$-algebra and differential field generated by $a$, respectively.
A family $(a_j)_{j\in I}$ is \emph{differentially algebraically
independent} over $K$ if the family  
$(\delta^\nu a_j)_{j\in I,\,\nu\in\NN^m}$ is algebraically
independent over $K$. The \emph{differential transcendence degree}
$\operatorname{dtrdeg}_K L$ of a differential field extension
$L/K$ is the cardinality of a maximal differentially algebraically
independent subset of $L$ over $K$ \cite[Ch.~II, \S9]{kolchin}. 

 We shall use the following consequence of additivity of differential
transcendence degree. 

\begin{lemma}\label{lem:free}
Let $T$ be a differentially algebraically independent family over a
differential field $E$ of characteristic zero. Every element of
$E\langle T\rangle\setminus E$ is differentially transcendental over $E$.
\end{lemma}

\begin{proof}
 Suppose, towards a contradiction, that
$r\in E\langle T\rangle\setminus E$ is differentially algebraic
over $E$. Choose a finite subtuple $t=(t_1,\ldots,t_n)$ of $T$
and write
\[
    r=A(t)/B(t),
    \qquad A,B\in E\{X_1,\ldots,X_n\},\quad B(t)\ne0.
\]
The polynomial $A-rB$ is nonzero. Otherwise comparison at a
monomial with nonzero coefficient in $B$ would give $r\in E$.
Since $(A-rB)(t)=0$, the tuple $t$ is differentially dependent
over $E\langle r\rangle$. Additivity in the tower
$E\subseteq E\langle r\rangle\subseteq E\langle t\rangle$ gives
\[
 \begin{aligned}
 n=\operatorname{dtrdeg}_E E\langle t\rangle
 &=\operatorname{dtrdeg}_E E\langle r\rangle
   +\operatorname{dtrdeg}_{E\langle r\rangle} E\langle t\rangle 
  \le 0+(n-1),
 \end{aligned}
\]
a contradiction. 
\end{proof}

  An ideal is \emph{differential} if it is closed under every
$\delta_i$. We write $[S]$ for the differential ideal generated by
a subset $S$. A prime or radical differential ideal is a
differential ideal that is prime or radical as an ordinary ideal.
In a differential $K$-algebra, $\sqrt{[S]}$ is the radical
differential ideal generated by $S$. For an ideal $I$ and an
element $h$, its saturation is
\[
       I:h^\infty=\{Q:h^eQ\in I\text{ for some }e\ge0\}.
\]
A \emph{differential affine $K$-domain} is an integral domain
finitely generated as a differential $K$-algebra. 

The closed sets of the \emph{Kolchin topology} on $\U^n$
are
\[
    V_\Delta(S)
    =
    \{a\in\U^n:\,\,H(a)=0\text{ for every }H\in S\},
    \quad
    S\subseteq\U\{X_1,\ldots,X_n\}.
\]
Restricting the coefficients to $K$ gives the Kolchin topology
over $K$.  On $\U^n\times\PP^1(\U)$, closedness is tested on the two
standard affine charts. Their coordinate change $W=1/Z$ preserves
Kolchin closed sets on the overlap, since denominators can be cleared
there.

 For a tuple $x$ of length $\ell$ in a differential field extension
of $K$, write $I_K(x)$ for the kernel of evaluation at $x$ on
$K\{X_1,\ldots,X_\ell\}$. It is a prime differential ideal.
A tuple $y$ of the same length is a \emph{differential
specialization of $x$ over $K$}, written $x\rightsquigarrow_K y$,
if every differential polynomial over $K$ that vanishes at $x$
also vanishes at $y$. Equivalently,
\[
             x\rightsquigarrow_K y\quad\Longleftrightarrow\quad
             I_K(x)\subseteq I_K(y).
\] 

\begin{definition} \label{def:graph-ideal}
 Let $f=F/G$, where $F,G\in K\{Y\}$ are nonzero and relatively prime.
The quotient is  \emph{indeterminate} at a common zero of $F$ and $G$.
Choose a differentially algebraically independent tuple
$t=(t_1,\ldots,t_n)$ over $K$ and put $u=f(t)$.
The \emph{graph ideal} of $f$ in the affine chart with
coordinate $Z$ is
\[
    J=I_K(t,u)
    =\big\{P\in  K\{Y,Z\}: \, P(t,u)=0\big\}.
\]
In the reciprocal chart with coordinate $W$, the graph ideal is
\[
    J^\vee=I_K(t,1/u)\subseteq K\{Y,W\}.
\] 
\end{definition}

 The two graph ideals depend only on $f$. Inverting $G$ identifies
the quotient of $K\{Y,Z\}[1/G]$ by the differential ideal generated
by $GZ-F$ with $K\{Y\}[1/G]$ via $Z\mapsto F/G$. Contracting to
$K\{Y,Z\}$, and then interchanging $F$ and $G$ give the following two identities
\[
    J=[GZ-F]:G^\infty,\qquad J^\vee=[FW-G]:F^\infty.
\]

\begin{definition} \label{def:projective-fiber}
For $f=F/G$ as above, let
\[
    \Gamma_f
    =
    \{(a,[F(a):G(a)]):
      a\in\U^n,\ G(a)\ne0\}
    \subseteq\U^n\times\PP^1(\U).
\]
The \emph{projective fiber} of $f$ at $b\in\U^n$ is
the set
\[
    \operatorname{Fib}_b(f)
    :=
    \{p\in\PP^1(\U):
      (b,p)\in\overline{\Gamma_f}^{\,\mathrm{Kol}}\},
\]
where the closure is taken in the Kolchin topology. 
\end{definition}

  Write $[T_0:T_1]$ for the homogeneous coordinates on $\PP^1$.
In the finite chart $Z=T_0/T_1$, the set $\Gamma_f$ is the graph of
$Z=F/G$ over $G\ne0$. In the reciprocal chart $W=T_1/T_0$, the
function $W=G/F$ is defined on all of $F\ne0$. A point of
$\Gamma_f$ lies in this chart precisely when $F(a)\ne0$ as well as
$G(a)\ne0$, since Definition~\ref{def:projective-fiber} already
requires $G(a)\ne0$. Thus the intersection of $\Gamma_f$ with the
reciprocal chart is the graph of $G/F$ restricted to $FG\ne0$.
The full graph of $G/F$ also contains the points with $F(a)\ne0$,
$G(a)=0$ and $W=0$. 

 The two differential substitution maps fixing $Y$ are
\[
\begin{aligned}
 K\{Y,Z\}&\longrightarrow K\{Y\}[1/G],&\qquad Z&\longmapsto F/G,\\
 K\{Y,W\}&\longrightarrow K\{Y\}[1/F],& W&\longmapsto G/F.
\end{aligned}
\]
Their kernels are $J$ and $J^\vee$, respectively. Tensoring over
$K$ with $\U$ is exact and commutes with localization, so the
kernels over $\U$ are the ideals generated by $J$ and $J^\vee$ in
the respective differential polynomial rings. 

  To see why restricting the reciprocal graph to $FG\ne0$
leaves its equations unchanged, take $Q\in\U\{Y,W\}$ and write
$Q(Y,G/F)=H/F^e$ with $H\in\U\{Y\}$. If $H\ne0$, then $HFG$
is a nonzero differential polynomial. Its nonvanishing at a tuple
of differential indeterminates and the differential closedness of
$\U$ give $a\in\U^n$ with $H(a)F(a)G(a)\ne0$. At this point
$Q(a,G(a)/F(a))\ne0$. Consequently, if $Q$ vanishes on the graph
over $FG\ne0$, then $H=0$ and $Q$ vanishes on the graph over all
of $F\ne0$. The converse follows by restriction. The same argument
in the finite chart uses the denominator $G$. Since taking closure
commutes with restriction to an open chart, the two charts of
$\overline{\Gamma_f}^{\,\mathrm{Kol}}$ are $V_\Delta(J)$ and
$V_\Delta(J^\vee)$, respectively. 

 Identifying $z\in\U$ with $[z:1]$ and writing
$\infty=[1:0]$, we obtain
\begin{equation}\label{eq:fiber-charts}
\begin{aligned}
    [z:1]\in\operatorname{Fib}_b(f)
    &\quad\Longleftrightarrow\quad
    (t,u)\rightsquigarrow_K(b,z)
    \qquad(z\in\U),\\
    \infty\in\operatorname{Fib}_b(f)
    &\quad\Longleftrightarrow\quad
    (t,1/u)\rightsquigarrow_K(b,0).
\end{aligned}
\end{equation}
Thus the finite points of the fiber are the common zeros of the
specialized polynomials $P(b,Z)$ for $P\in J$. At a common zero
of $F$ and $G$, the equation $GT_0-FT_1=0$ holds for every
$[T_0:T_1]$, even when the graph ideals determine a unique point. 

\subsection{Valuation rings and residue fields}
\label{subsec:valuation-rings}

 An additive valuation on a field $E$ takes values in
$\Gamma\cup\{\infty\}$, where $\Gamma$ is a totally ordered
abelian group. It satisfies $v(x)=\infty$ exactly when $x=0$, and
\[
    v(xy)=v(x)+v(y),
    \qquad
    v(x+y)\geq\min\{v(x),v(y)\}.
\]
We order $\infty$ above $\Gamma$ and set $\infty+\alpha=\infty$ for
$\alpha\in\Gamma\cup\{\infty\}$, and take
$\Gamma=v(E^\times)$ as the \emph{value group}.
The valuation is \emph{real} if its value group is an
ordered subgroup of $\RR$.
The valuation ring, maximal ideal, and residue field of $(E, v)$ are
\[
    \OO_v=\{x\in E:v(x)\geq0\},
    \qquad
    \mm_v=\{x\in E:v(x)>0\},
    \qquad
    k_v=\OO_v/\mm_v.
\]
We write $\res_v:\OO_v\to k_v$ for the residue map.
An element of $\OO_v$ is a unit precisely when its value is zero. For a subring $B\subseteq\OO_v$, the
\emph{center} of $v$ on $B$ is the prime ideal
$\mm_v\cap B$. 

 A \emph{semivaluation} on a ring $R$ satisfies the same product
and sum rules, with $v(0)=\infty$ and $v(1)=0$, but its
\emph{support}
\(
    \{x\in R:v(x)=\infty\}
\)
may contain nonzero elements. For valuations and
semivaluations alike,
\[
    v(x)<v(y)\quad\Longrightarrow\quad v(x+y)=v(x).
\]
Consequently, a finite sum with a unique summand of least value
has the value of that summand. A semivaluation with zero support
on a domain extends uniquely to a valuation on its fraction field.
We also call its restriction to the domain a valuation. 

A valuation on a differential field $(E,\Delta)$ is
\emph{monotone} if
\[
    v(\delta_i x)\geq v(x)
    \qquad(x\in E,\ 1\leq i\leq m).
\]
Then $\OO_v$ is a differential subring and $\mm_v$ is a
differential ideal. The derivations descend to $k_v$,
and $\res_v$ is a differential homomorphism. We use the
same monotonicity condition for semivaluations on
differential rings.

A local ring $(V,\mm_V)$ \emph{dominates} a local ring
$(R,\mm_R)$ if $R\subseteq V$ and
$\mm_V\cap R=\mm_R$. A subring $V$ of a field $E$ is a
\emph{valuation ring} if, for every $x\in E^\times$,
at least one of $x$ and $x^{-1}$ belongs to $V$.
Its associated valuation has value group
$E^\times/V^\times$, written additively and ordered by
\[
    v(x)\geq v(y)
    \quad\Longleftrightarrow\quad
    x/y\in V
    \qquad(x,y\in E^\times).
\]
In particular, $V=\OO_v$.

The construction in Section~\ref{sec:valuation} relies on Chevalley's extension theorem and H\"older's embedding theorem.

\begin{classicaltheorem}[Chevalley's extension theorem {\cite[Thm.~3.1.1]{engler-prestel}}]
  Let $R$ be a subring of a field $E$, and let $\mathfrak r$ be a
prime ideal of $R$. There exists a valuation ring $V$ of $E$
dominating $R_{\mathfrak r}$. In particular, $R\subseteq V$ and
$\mm_V\cap R=\mathfrak r$. 
\end{classicaltheorem}

 To pass from an arbitrary value group to real-valued valuations, we use convex subgroups. Recall that a subgroup $I$ of a totally ordered abelian group $\Gamma$ is \emph{convex} if $0\le\alpha\le\beta$ and $\beta\in I$ imply $\alpha\in I$. The quotient $\Gamma/I$ carries the induced total order given by
\[
    \alpha + I > 0 + I
    \quad\Longleftrightarrow\quad
    \alpha > 0 \ \text{and}\ \alpha \notin I.
\]
A totally ordered abelian group is \emph{archimedean} if for every pair of elements $\alpha, \beta > 0$, there is a positive integer $r$ with $r\alpha \geq \beta$. 

\begin{classicaltheorem}[H\"older's theorem {\cite[Prop.~2.1.1]{engler-prestel}}]
Let $\Gamma$ be an archimedean ordered abelian group. Given $\gamma\in\Gamma_{>0}$, there is an order-preserving group embedding $\iota:\Gamma\hookrightarrow\mathbb R$ with $\iota(\gamma)=1$.
\end{classicaltheorem}

\section{A valuative generization of a differential prime}\label{sec:valuation}

In this section, we shall construct a monotone semivaluation on a differential affine $K$-domain
whose positive-value ideal is a prescribed nonzero prime differential
ideal. We first normalize a chosen positive value in an arbitrary
ordered value group. 

\begin{lemma}\label{lem:truncation}
  Let $v$ be a valuation with valuation ring $V$ and value group
$\Gamma$, and let $\gamma\in\Gamma_{>0}$. There exist a prime ideal
$Q$ of $V$ and a nonnegative real valuation $\lambda$ on $V/Q$
with the following properties for $x\in V$. 
\begin{enumerate}
  \item[(i)]
If $v(x)=0$, then $\lambda(x+Q)=0$.
  \item[(ii)]
If $v(x)\ge\gamma$, then $\lambda(x+Q)\ge1$.
  \item[(iii)]
If $v(x)=\gamma$, then $\lambda(x+Q)=1$.
\end{enumerate}
In particular, $x\notin Q$ in cases \textup{(i)} and \textup{(iii)}.
\end{lemma}

\begin{proof}
  For $\beta\in\Gamma$ write $|\beta|=\max\{\beta,-\beta\}$. Set
\[
 \begin{aligned}
 H&=\{\beta:|\beta|\le r\gamma\text{ for some integer }r\ge1\},\\
 I&=\{\beta:r|\beta|<\gamma\text{ for every integer }r\ge1\}.
 \end{aligned}
\]
Both sets are convex and stable under negation. Closure of $H$ under
addition follows from the triangle inequality in $\Gamma$. For
$\alpha,\beta\in I$ and $r\ge1$ we have
$r|\alpha+\beta|\le2r\max\{|\alpha|,|\beta|\}<\gamma$, proving
closure of $I$ under addition. Thus $I\subseteq H$ are convex subgroups
and $\gamma\in H\setminus I$. 

  To verify that $H/I$ is archimedean, take positive cosets represented
by $\alpha,\beta>0$ in $H\setminus I$. Since $\alpha\notin I$,
some multiple $r\alpha$ is at least $\gamma$. Since $\beta\in H$,
some multiple $s\gamma$ is at least $\beta$. Thus
$\beta\le s\gamma\le sr\alpha$, and the same inequality holds in
the quotient. By H\"older's theorem there is an order-preserving
embedding $\iota:H/I\hookrightarrow\RR$ sending $\gamma+I$ to $1$. 

  Let $Q=\{x\in V:v(x)>h\text{ for every }h\in H\}$. The valuation
inequality and nonnegativity on $V$ make $Q$ a proper ideal. For
$x\in V\setminus Q$, convexity gives $v(x)\in H$, since
$0\le v(x)\le h$ for some $h\in H$. Hence $V\setminus Q$ is
multiplicatively closed and $Q$ is prime. On $V/Q$ define
$\lambda(0)=\infty$ and
\[
  \lambda(x+Q)=\iota(v(x)+I)\qquad(x\notin Q).
\]
Adding an element of $Q$ to $x\notin Q$ preserves $v(x)$, so
$\lambda$ is well defined. The valuation laws follow from those for
$v$ and the order-preserving homomorphism $H\to\RR$ induced by
$\iota$. The values $0$ and $\gamma$ map to $0$ and $1$,
respectively, which proves \textup{(i)} and \textup{(iii)}. If
$v(x)\ge\gamma$, then either $x\in Q$, giving $\lambda(x+Q)=\infty$,
or order preservation gives $\lambda(x+Q)\ge1$. 
\end{proof}

 We next give the variables of a formal power series a common
positive weight. This makes the minimum in the Gauss valuation
attainable and bounds the loss of value under differentiation. 

\begin{lemma}\label{lem:gauss}
Let $D$ be a domain with a nonnegative real valuation $\lambda$
such that $\lambda(n\cdot 1)=0$ for every nonzero integer $n$.
For $\varepsilon>0$ put $W_\varepsilon(0)=\infty$ and, on nonzero
series in $D[[X_1,\ldots,X_m]]$, set
\[
 W_\varepsilon\!\left(\sum_\nu b_\nu X^\nu\right)
          =\min_\nu\{\lambda(b_\nu)+\varepsilon|\nu|\}.
\]
Then $W_\varepsilon$ is a nonnegative real valuation, and
\[
 W_\varepsilon\!\left(\frac{\partial f}{\partial X_i}\right)
    \ge W_\varepsilon(f)-\varepsilon
    \qquad (1\le i\le m).
\]
\end{lemma}

\begin{proof}
 For a nonzero series $f=\sum_\nu b_\nu X^\nu$, let $M$ be the
weight of any nonzero term. Since $\lambda$ is nonnegative, a term
of weight at most $M$ must have $|\nu|\le M/\varepsilon$. There are
only finitely many such indices, so the minimum is attained. The
inequality for sums follows coefficient by coefficient. 

  Let $g=\sum_\mu c_\mu X^\mu$ be another nonzero series.  Every product term $b_\nu c_\mu X^{\nu+\mu}$ has weight
\(
 \lambda(b_\nu c_\mu)+\varepsilon|\nu+\mu|
 \ge W_\varepsilon(f)+W_\varepsilon(g),
\)
so $W_\varepsilon(fg)\ge W_\varepsilon(f)+W_\varepsilon(g)$. Let $\alpha$ and $\beta$ be the lexicographically smallest multi-indices attaining $W_\varepsilon(f)$ and $W_\varepsilon(g)$, respectively. In the coefficient $\sum_{\nu+\mu=\alpha+\beta} b_\nu c_\mu$ of $X^{\alpha+\beta}$, the term $b_\alpha c_\beta$ satisfies
\[
 \lambda(b_\alpha c_\beta)= \lambda(b_\alpha)+\lambda(c_\beta)=W_\varepsilon(f)+W_\varepsilon(g)-\varepsilon|\alpha+\beta|.
\]
For any other pair $(\nu,\mu)$ with $\nu+\mu=\alpha+\beta$, either $\nu<_{\mathrm{lex}}\alpha$ or $\mu<_{\mathrm{lex}}\beta$, so $\lambda(b_\nu c_\mu)+\varepsilon|\alpha+\beta|>W_\varepsilon(f)+W_\varepsilon(g)$, whence $\lambda(b_\nu c_\mu)>\lambda(b_\alpha c_\beta)$. By the strict triangle inequality for $\lambda$, the coefficient of $X^{\alpha+\beta}$ has $\lambda$-value $\lambda(b_\alpha c_\beta)$, giving $W_\varepsilon(fg)\le W_\varepsilon(f)+W_\varepsilon(g)$ and hence $W_\varepsilon(fg)=W_\varepsilon(f)+W_\varepsilon(g)$.

 Finally, a term with $\nu_i>0$ differentiates to
$\nu_i b_\nu X^{\nu-e_i}$. Since $\lambda(\nu_i\cdot1)=0$, its
weight decreases by exactly $\varepsilon$. This proves the derivative
bound.
\end{proof}

We now construct a semivaluation on $A$ whose positive-value ideal
is the prescribed prime $\pp$. A uniform bound on a finite set of
generators will force its support to be strictly smaller than $\pp$.

\begin{theorem}\label{lem:valuation}
Let $A$ be a differential affine $K$-domain and let $\pp$ be a
nonzero prime differential ideal of $A$. There exist a prime
differential ideal $\qq\subsetneq\pp$ and a monotone real valuation
$v$ on $L=\Frac(A/\qq)$ such that $A/\qq\subseteq\OO_v$ and
\[
                      \mm_v\cap(A/\qq)=\pp/\qq .
\]
In particular, the composite differential homomorphism
\[
 A\longrightarrow A/\qq\hookrightarrow\OO_v
              \xrightarrow{\ \res_v\ } k_v
\]
has kernel $\pp$.
\end{theorem}

\begin{proof}
  By the Ritt--Raudenbush basis theorem \cite[Ch.~III]{kolchin}, choose
nonzero elements $h_1,\ldots,h_\ell$ generating $\pp$ as a radical
differential ideal. We shall construct a nonnegative monotone
semivaluation $w$ on $A$ satisfying
\[
 w(a)=0\quad(a\notin\pp),\qquad 1\le\min_j w(h_j)\le2.
\]
The lower bound will determine its positive-value ideal, and the
upper bound will ensure that some $h_j$ lies outside its support. 

For each $N\ge1$ we construct a semivaluation $w_N$
for which differentiation lowers values by at most $1/N$.
By Chevalley's extension theorem, there is a valuation ring $V$ of $\Frac(A)$
containing $A$ with $\mm_V\cap A=\pp$. Its valuation $v_0$ is
positive on $\pp\setminus\{0\}$ and zero on $A\setminus\pp$, in
particular on $K^\times$. Put
\[
 \gamma_N=\min\{v_0(\delta^\nu h_j):1\le j\le\ell,\ |\nu|\le N,\
                                     \delta^\nu h_j\ne0\}.
\]
This is the minimum of a finite nonempty set of positive values,
since $\pp$ is differential.
Apply Lemma~\ref{lem:truncation} with $\gamma=\gamma_N$. Write
$D_N=V/Q_N$ for the resulting domain, $\lambda_N$ for its real
valuation and $\pi_N:A\to D_N$ for the quotient map restricted to
$A$. By part~\textup{(i)} of that lemma, $\pi_N$ embeds $K$ in
$D_N$, and every nonzero element of $K$ has value zero there.

The Taylor homomorphism associated to $\pi_N$ \cite[Example~4.3]{ng}
is given by
\begin{equation}\label{eq:taylor}
 T_N(a)=\sum_{\nu\in\NN^m}
       \pi_N\!\left(\frac{\delta^\nu a}{\nu!}\right) X^\nu
       \quad\text{in }D_N[[X_1,\ldots,X_m]].
\end{equation}
The Leibniz formula gives multiplicativity, and termwise
differentiation gives $T_N(\delta_i a)=\partial T_N(a)/\partial X_i$.
Composing $T_N$ with the weighted Gauss valuation of
Lemma~\ref{lem:gauss} for $\varepsilon=1/N$ gives a nonnegative
semivaluation $w_N=W_{1/N}\circ T_N$ on $A$ satisfying
\begin{equation}\label{eq-diffbound}
 w_N(\delta_i a)\ge w_N(a)-1/N\qquad(a\in A,\ 1\le i\le m).
\end{equation}
For $a\notin\pp$,  $\pi_N(a)$ has value
zero by Lemma~\ref{lem:truncation}\textup{(i)}, so $w_N(a)=0$.

For each $j$ and $|\nu|\le N$, the coefficient indexed by $\nu$ in
$T_N(h_j)$ has $\lambda_N$-value at least $1$. Indeed, whenever
$\delta^\nu h_j\ne0$ we have
\[
 v_0(\delta^\nu h_j/\nu!)=v_0(\delta^\nu h_j)\ge\gamma_N,
\]
since $\nu!\in K^\times$ has $v_0$-value zero. A zero coefficient
has value $\infty$. For the remaining coefficients the bound follows from
Lemma~\ref{lem:truncation}\textup{(ii)}. Terms of higher degree have
weight at least $|\nu|/N>1$. Thus $w_N(h_j)\ge1$ for every $j$.
For some $j_0,\nu_0$ the value $v_0(\delta^{\nu_0}h_{j_0})$ equals
$\gamma_N$ with $|\nu_0|\le N$. By
Lemma~\ref{lem:truncation}\textup{(iii)}, the corresponding term
has weight $1+|\nu_0|/N\le2$. Hence
\begin{equation}\label{eq:bounds}
 1\le\min_{1\le j\le\ell} w_N(h_j)\le2 .
\end{equation}

To remove the error $1/N$ in~\eqref{eq-diffbound} while
preserving the bounds in~\eqref{eq:bounds}, we use compactness.
Set $u_N(a)=e^{-w_N(a)}$, with $e^{-\infty}=0$. Regard $u_N$ as a
point of $X=[0,1]^A$, the space of functions $A\to[0,1]$ with the
product topology. By Tychonoff's theorem $X$ is compact. The closed
sets
\[
 C_M=\overline{\{u_N:N\ge M\}}\qquad(M\ge1)
\]
are nonempty and nested, so every finite subfamily has nonempty
intersection. Compactness therefore gives $u\in\bigcap_{M\ge1}C_M$.

Every closed condition satisfied by all $u_N$ for $N\ge M$ holds on
$C_M$ and hence at $u$. For example, for fixed $a,b\in A$, the
condition $f(ab)=f(a)f(b)$ is closed because the coordinate maps and
multiplication are continuous. The same argument applies to the sum
inequality, the prescribed values, and the finite maximum in
\eqref{eq:bounds}. Consequently
\[
 \begin{gathered}
 u(0)=0,\quad u(1)=1,\quad u(ab)=u(a)u(b),\quad
 u(a+b)\le\max\{u(a),u(b)\},\\
 u(a)=1\quad(a\notin\pp),\qquad
 e^{-2}\le\max_j u(h_j)\le e^{-1}.
 \end{gathered}
\]
Fix $a\in A$, $i$ and $M\ge1$. By \eqref{eq-diffbound},
$u_N(\delta_i a)\le e^{1/M}u_N(a)$ for all $N\ge M$. The inequality
$f(\delta_i a)\le e^{1/M}f(a)$ is again a closed condition, so
$u(\delta_i a)\le e^{1/M}u(a)$. Letting $M\to\infty$ yields
$u(\delta_i a)\le u(a)$. Thus $w=-\log u$, with $-\log 0=\infty$,
has the required monotonicity, semivaluation laws and bounds.

 Let $\qq=\{a\in A:w(a)=\infty\}$ be the support of $w$,
and put $I=\{a\in A:w(a)>0\}$. Nonnegativity
and the semivaluation laws make both sets prime ideals, and monotonicity
makes them differential. Since $I$ contains all the $h_j$, it
contains $\pp$. Since $w$ vanishes off $\pp$, we also have
$I\subseteq\pp$, whence $I=\pp$. The upper bound on $\min_j w(h_j)$
gives some $h_j\in\pp\setminus\qq$, proving $\qq\subsetneq\pp$.

Adding an element of $\qq$ to an element of finite value preserves
its value. Thus $w$ induces a valuation on $A/\qq$, finite on every
nonzero element, and extends by division to a real valuation $v$ of
$L=\Frac(A/\qq)$. For nonzero $x,y\in A/\qq$, monotonicity on
$A/\qq$ gives
\[
 \frac{\delta_i(x/y)}{x/y}
    =\frac{\delta_i x}{x}-\frac{\delta_i y}{y}\in\OO_v.
\]
Hence $v(\delta_i f)\ge v(f)$ for every $f\in L$. Nonnegativity on
$A/\qq$ gives $A/\qq\subseteq\OO_v$, and the residue of $a+\qq$ is
zero exactly when $w(a)>0$, that is, when $a\in\pp$. This proves
the center and kernel assertions
\end{proof}

\begin{remark} 
\label{rem:ordinary-center}
For one derivation, Fang, Freitag and Harrison-Trainor
use monotone valuations to study containment of prime differential
ideals. Their construction, which draws on Morrison's work on
neighboring differential primes \cite[\S2]{morrison}, also yields
the prescribed center in Theorem~\ref{lem:valuation}
\cite[proof of Lemma~6.2]{ffh}. The Taylor construction above gives
a single valuation monotone for every derivation in $\Delta$.
\end{remark}

\section{Lowest-degree terms under valued specialization}\label{sec:comparison}

  Let $J=I_K(t,u)$ be the graph ideal of Definition~\ref{def:graph-ideal},
where $t$ is differentially algebraically independent over $K$ and
$u=F(t)/G(t)$. Every $P\in J$ satisfies $P(Y,F/G)=0$ in
$K\{Y\}[1/G]$. Hence $P(y,F(y)/G(y))=0$ at every tuple $y$ in a
differential field extension of $K$ for which $G(y)\ne0$. We use
these identities to compare the lowest homogeneous terms of $(F,G)$
at a point and at its residue. 

  For a tuple $T=(T_1,\ldots,T_n)$ of differential indeterminates,
expand the translated pair in differential monomials of $T$ as  
\begin{equation}\label{eq:expansion}
 (F(Y+T),G(Y+T))=\sum_M (f_M(Y),g_M(Y))\,M(T).
\end{equation}
The sum runs over finitely many distinct differential monomials in $T$, including those of order zero. Here $f_M,g_M\in K\{Y\}$ and
$d_M$ is the total degree of $M$, with each derivative assigned
degree one. At a point $a\in E^n$ in a differential field extension
$E/K$, the coefficient pair of $M(T)$ is $(f_M(a),g_M(a))$.
Translation by $a$ is an automorphism of $E\{T\}$ with inverse
translation by $-a$, so $F(a+T)$ and $G(a+T)$ are both nonzero.
Set
\[
 d(a)=\min\{d_M:(f_M(a),g_M(a))\ne(0,0)\}.
\]
Write $F_{a,d}$ and $G_{a,d}$ for the homogeneous parts of degree
$d=d(a)$ in $F(a+T)$ and $G(a+T)$, respectively. At least one is
nonzero. If $s$ is an indeterminate with $\delta_i s=0$ for all $i$,
then $M(sT)=s^{d_M}M(T)$, so these parts are the coefficients of
$s^d$ in $F(a+sT)$ and $G(a+sT)$. 

\begin{definition}
 The pair of lowest homogeneous parts defines a projective point
\[
            \rho(a;T)=[F_{a,d}:G_{a,d}]\in\PP^1(E\langle T\rangle).
\]
If this point is defined over $E$, we call it the \emph{radial value}
at $a$. Explicitly, for $[p:q]\in\PP^1(E)$ the equality
$\rho(a;T)=[p:q]$ means that $qF_{a,d}=pG_{a,d}$. Comparing
coefficients gives
\begin{equation}\label{eq-rho}
 \rho(a;T)=[p:q]\,\,\Longleftrightarrow\,\,
 qf_M(a)=pg_M(a)\,\,\text{for every }M\text{ with }d_M=d(a).
\end{equation}
Thus $\rho(a;T)$ is defined over $E$ precisely when the nonzero
coefficient pairs of degree $d(a)$ span a single line in $E^2$. 
\end{definition}

 Let $(L,v)$ be a valued field with residue field $k$, and write
$\bar x$ for the residue of $x\in\OO_v$. For $[A:B]\in\PP^1(L)$,
choose an entry $e\in\{A,B\}$ of minimum value. Then $A/e$ and $B/e$
lie in $\OO_v$, and one equals $1$. Define the projective reduction by
\begin{equation}\label{eq-res-reduction}
           \res_v[A:B]=[\overline{A/e}:\overline{B/e}]
                    \in\PP^1(k).
\end{equation}
If both entries have minimum value, the two normalized pairs differ
by a unit and hence represent the same residue point. Multiplying
$A,B,e$ by a common nonzero scalar leaves the normalized pair
unchanged. Thus \eqref{eq-res-reduction} depends only on $[A:B]$. 

  For a nonzero coefficient pair $(A_M,B_M)$, we call
$\res_v[A_M:B_M]$ its \emph{reduced direction}. To compare the
directions arising in \eqref{eq:expansion}, we perturb $a$ by $sT$
and vary the weight $q>0$ assigned to $s$. In the Gauss valuation
$v_q$ this means $v_q(s)=q$. The following lemma provides a valued
extension in which the derivatives of $T$ remain algebraically
independent after reduction. 

\begin{lemma}\label{lem:gauss-ext}
  Let $(L,v)$ be a monotone real valued differential field with residue
field $k$, and fix $q>0$. Let $s$ be a constant indeterminate over
$L$, and let $T=(T_1,\ldots,T_n)$ be differential indeterminates over
$L(s)$. There is a monotone real valuation $v_q$ on
$L(s)\langle T\rangle$ extending $v$ with the following properties. 
\begin{enumerate}
\item[\textup{(i)}]
 For every nonzero polynomial $\sum_j a_js^j\in L[s]$,
\[
    v_q\!\left(\sum_j a_js^j\right)
    =\min_{j:a_j\ne0}\bigl(v(a_j)+jq\bigr).
\] 
\item[\textup{(ii)}]
  The value of a nonzero $H\in L(s)\{T\}$ is the minimum of the
$v_q$-values of its nonzero coefficients. 
\item[\textup{(iii)}]
  Write $k_q$ for the residue field of $(L(s),v_q)$. There is a
natural differential embedding $k\hookrightarrow k_q$. The tuple
$\bar T$ is differentially algebraically independent over $k_q$, and
the residue field of $(L(s)\langle T\rangle,v_q)$ is
$k_q\langle\bar T\rangle$. 
\end{enumerate}
\end{lemma}

\begin{proof}
 Give $s$ weight $q$ in $L[s]$. Use the minimum in
\textup{(i)} to assign values to polynomials of $L[s]$ and extend by division
to $L(s)$. Over this field give each jet $\delta^\nu T_j$ weight
zero and use the minimum in \textup{(ii)}. To verify the valuation
laws at both stages, observe that a nonzero polynomial has finite
support, so the minimum term weight exists for arbitrary coefficient
values. The sum inequality holds coefficientwise. In each of two
nonzero factors, choose the lexicographically first exponent of
minimum weight, say $\alpha$ and $\beta$, ordering the finitely
many variables involved. Among the products contributing to the
coefficient of exponent $\alpha+\beta$, only the chosen product
can have minimum weight. For any other decomposition
$\nu+\mu=\alpha+\beta$, either $\nu$ precedes $\alpha$ or $\mu$
precedes $\beta$, so one factor has strictly greater weight.
The monomial weight is the same for all terms contributing to this
coefficient. Its unique summand of least value therefore survives
addition, proving multiplicativity. Both extensions by division
are thus well defined, the second giving $v_q$ on
$L(s)\langle T\rangle$.

On $L[s]$, differentiation acts on the coefficients and annihilates $s$,
so monotonicity follows from that of $v$. The quotient rule extends
the inequality to $L(s)$. In $L(s)\{T\}$ each derivative of a
monomial is a sum of monomials of weight zero with integer
coefficients of nonnegative value. Thus differentiation cannot lower
the least coefficient value. Another application of the quotient
rule gives monotonicity on $L(s)\langle T\rangle$. 

  The residue maps therefore commute with every derivation, and the
restriction of $v_q$ to $L$ induces the differential embedding
$k\hookrightarrow k_q$. For a nonzero differential polynomial over
$k_q$ in $n$ indeterminates, lift its nonzero coefficients to units
of the valuation ring of $L(s)$ and evaluate the lift at $T$. By
\textup{(ii)} the result has value zero, so its residue, the value
of the original polynomial at $\bar T$, is nonzero. This proves
differential algebraic independence of $\bar T$. 

 Take an integral fraction $H_1/H_2$ with $H_1,H_2\in L(s)\{T\}$.
Choose a coefficient $e$ of least value in its denominator. Then
$v_q(H_1)\ge v_q(H_2)=v_q(e)$. Dividing both polynomials by $e$
gives integral coefficients and makes the denominator a unit. Their
reductions therefore give
\[
 \res_{v_q}(H_1/H_2)=\frac{\overline{H_1/e}}{\overline{H_2/e}}
    \in k_q\langle\bar T\rangle .
\]
Conversely, the residue field contains $k_q$ and all
$\delta^\nu\bar T_j=\overline{\delta^\nu T_j}$, so it contains the
field they generate. This proves \textup{(iii)}. 
\end{proof}

  If the smaller value of the entries in a coefficient pair is $c_M$,
then the substitution $T\mapsto sT$ changes that value to
$c_M+d_Mq$. Different indices can attain the minimum as $q$ varies.
The next lemma compares their labels. 

\begin{lemma}\label{lem:envelope}
Let $\mathcal M$ be a finite nonempty index set. To each $M\in\mathcal M$
assign $c_M\in\RR$, $d_M\in\NN$ and a label $D_M$ in a set $X$.
For $q>0$ write $\ell_M(q)=c_M+d_Mq$ and put
\[
 \mathcal A(q)=
 \{M\in\mathcal M:\ell_M(q)=\min_{N\in\mathcal M}\ell_N(q)\}.
\]
Suppose that, for each $q>0$, the labels of the indices in
$\mathcal A(q)$ coincide. Then there is a single $D\in X$ such
that $D_M=D$ for all $M\in\bigcup_{q>0}\mathcal A(q)$.
\end{lemma}

\begin{proof}
 Let $D(q)$ denote the common label of the indices in
$\mathcal A(q)$. Fix $q_0>0$ and choose $M_0\in\mathcal A(q_0)$.
For each $M\notin\mathcal A(q_0)$, the continuous function
$g_M(q)=\ell_M(q)-\ell_{M_0}(q)$ is positive at $q_0$.
Since there are only finitely many such $M$, there is an
$\varepsilon\in(0,q_0)$ for which all these differences remain
positive whenever $|q-q_0|<\varepsilon$. At each such $q$, every
index outside $\mathcal A(q_0)$ has greater weight than $M_0$,
so $\mathcal A(q)\subseteq\mathcal A(q_0)$. Choosing any
$M\in\mathcal A(q)$ now gives $D(q)=D_M=D(q_0)$. 

  Thus $D(q)$ is locally constant on $(0,\infty)$.
Let $D=D(1)$. The set $\{q>0:D(q)=D\}$ and its complement are
both open by local constancy. Since this set contains $1$ and
$(0,\infty)$ is connected, it is the whole interval. Hence $D_M=D$
for every index that attains the minimum at some $q>0$.
\end{proof}

  We apply these two lemmas to the coefficient pairs of the translated
polynomials at $a$ and at its residue $b$. A graph equation that
remains nonzero as a differential polynomial in the ordinate at $b$
will force the minimizing pairs to have a common reduced direction. 

\begin{proposition}\label{lem:comparison}
 Let $(L,v)$ be a monotone real valued differential field over $K$
with residue field $k$, and suppose $v(c)=0$ for every $c\in K^\times$. Let
$a\in\OO_v^n$ have residue $b\in k^n$. Suppose $P(b,Z)\ne0$ in
$k\{Z\}$ for some $P\in J$. If $\rho(a;T)\in\PP^1(L)$, then
$\rho(b;T)\in\PP^1(k)$ and
\[
                 \rho(b;T)=\res_v\bigl(\rho(a;T)\bigr).
\]
\end{proposition}

\begin{proof}
 Write $(A_M,B_M)=(f_M(a),g_M(a))$ in \eqref{eq:expansion}. Both
entries lie in $\OO_v$, since $a\in\OO_v^n$, $v$ is trivial on $K$,
and the valuation ring is stable under the derivations. Restrict the
indices to the nonzero pairs and put
\[
 c_M=\min\{v(A_M),v(B_M)\},\qquad D_M=\res_v[A_M:B_M]\in\PP^1(k).
\]
Each $c_M$ is nonnegative, and each $D_M$ is defined in the
original residue field $k$ before a weight is chosen. 

 Fix $q>0$. Let $v_q:\, L(s)\langle T\rangle\,\longrightarrow\, \mathbb R\cup\{\infty\}$ be the valuation extension of $v$ obtained by applying Lemma~\ref{lem:gauss-ext} to adjoin a differential
constant $s$ of value $q$ and a tuple $T$ of differential
indeterminates whose derivatives have value zero.  Write $k_q$ for the
residue field of its restriction to $L(s)$. The tuple $y=a+sT$
has residue $b$. Since $s$ is a differential
constant, the coefficient pair of $M(T)$ in $(F(y),G(y))$ is
$s^{d_M}(A_M,B_M)$. Its \emph{weight at $q$} is the smaller of
the valuations of its two entries. Since $v_q(s)=q$, it is given by
\[
 \ell_M(q)=\min\{v_q(s^{d_M}A_M),v_q(s^{d_M}B_M)\}=c_M+d_Mq .
\]
Let $m(q)=\min_M\ell_M(q)$. An index  $M$ is \emph{active at $q$} if
its pair attains this minimum, and we write
$\mathcal A(q)=\{M:\ell_M(q)=m(q)\}$. 

  Choose a coefficient $h$ of $F(y)$ or $G(y)$ with $v_q(h)=m(q)$.
Both polynomials have integral coefficients after division by $h$.
Reduce their coefficients and write the resulting differential
polynomials as $U(S),V(S)\in k_q\{S\}$, where $S$ is a tuple of
formal differential indeterminates. Thus $U(\bar T),V(\bar T)$
are the residues of $F(y)/h,G(y)/h$, respectively. Write $(u_M,v_M)$ for the
resulting coefficient pair at $M$. For an inactive index both
$s^{d_M}A_M/h$ and $s^{d_M}B_M/h$ have positive value, so both
coefficients reduce to zero. For an active index take an entry $e_M$
of $(A_M,B_M)$ with $v(e_M)=c_M$. Then
\[
 (u_M,v_M)=\overline{\frac{s^{d_M}e_M}{h}}\,
    \left(\overline{\frac{A_M}{e_M}},\overline{\frac{B_M}{e_M}}\right).
\]
The scalar is nonzero because $d_Mq+c_M=v_q(h)$, and one entry in
parentheses equals $1$. Thus the nonzero coefficient pairs of $U,V$
are precisely those indexed by active monomials, each with direction
$D_M$ in $\PP^1(k_q)$. 

  We use the graph equation to show that the pairs surviving
at $q$ have the same direction. If $V=0$, each
active pair represents $[1:0]$. Suppose $V\ne0$. By
Lemma~\ref{lem:gauss-ext}\textup{(iii)}, the tuple $\bar T$ is
differentially algebraically independent over $k_q$. It follows that
$\overline{G(y)/h}=V(\bar T)\ne0$. Hence $G(y)/h$ is a unit of
$\OO_{v_q}$. Since $F(y)/h$ also belongs to $\OO_{v_q}$, we have
\[
 \frac{F(y)}{G(y)}=\frac{F(y)/h}{G(y)/h}\in\OO_{v_q},\qquad
 r=\overline{\frac{F(y)}{G(y)}}=\frac{U(\bar T)}{V(\bar T)}
       \in k_q\langle\bar T\rangle .
\]
Reducing the graph identity $P(y,F(y)/G(y))=0$ gives $P(b,r)=0$.
As $P(b,Z)$ is nonzero in $k\{Z\}$ and $k$ embeds in $k_q$, the
element $r$ is differentially algebraic over $k_q$.
Lemma~\ref{lem:free} now gives $r\in k_q$. Thus $U-rV$ is a
differential polynomial over $k_q$ vanishing at $\bar T$, so it is
identically zero. Comparing coefficients gives $u_M=rv_M$. For an
active pair this forces $v_M\ne0$, and its projective class is
$[r:1]$. 

 In either case the points $D_M$ for $M\in\mathcal A(q)$
have the same image in $\PP^1(k_q)$. Since $k\hookrightarrow k_q$
is injective, they coincide in $\PP^1(k)$. Denote their common
direction by $D(q)$. 

 To see why $D(q)$ is independent of $q$, fix $q_0>0$.
If $M\notin\mathcal A(q_0)$, the continuous function
$\ell_M(q)-m(q)$ is positive at $q_0$ and remains positive
near $q_0$. There are only finitely many
indices, so $\mathcal A(q)\subseteq\mathcal A(q_0)$ throughout
some neighborhood of $q_0$, as in Lemma~\ref{lem:envelope}.
For each such $q$, choose $M\in\mathcal A(q)$. This pair is
active at both $q$ and $q_0$, and its direction $D_M$ gives
\[
                         D(q)=D_M=D(q_0).
\]
Thus $D(q)$ is locally constant on the connected interval
$(0,\infty)$. Write $D$ for this common direction. 

 To identify $D$ at $a$, take $q>\max_M c_M$. Every pair of
degree $d(a)$ then has weight less than $(d(a)+1)q$, and every pair
of higher degree has weight at least $(d(a)+1)q$. Thus each active
pair has degree $d(a)$. Since $\rho(a;T)$ is defined over $L$,
\eqref{eq-rho} shows that each such pair represents $\rho(a;T)$.
Reducing any one of them gives $D=\res_v(\rho(a;T))$. 

  At $b$, the coefficient pair $(f_M(b),g_M(b))=(\bar A_M,\bar B_M)$
is nonzero precisely when $c_M=0$. At least one pair survives,
since $G(T)\ne0$ in $k\{T\}$ and translation by $b$ is invertible.
By definition $d(b)$ is the least degree among the surviving pairs.
Their weights are $d_Mq$, with minimum $d(b)q$. Choose $q>0$ small
enough that $d(b)q<c_M$ for every positive $c_M$. This is possible
because the set of indices is finite. Each such index has weight
$c_M+d_Mq\ge c_M>d(b)q$. The active indices are therefore exactly
those of the surviving pairs of degree $d(b)$. For these indices,
$[f_M(b):g_M(b)]=[\bar A_M:\bar B_M]=D_M=D$. Applying
\eqref{eq-rho} at $b$ gives $\rho(b;T)=D=\res_v(\rho(a;T))$. 
\end{proof}

\begin{remark}\label{rem:change-degree}
Reduction can annihilate every coefficient of the lowest homogeneous
pair at $a$. The radial value at its residue $b$ is then determined
by terms of higher degree. For example, take $F=X^2$ and $G=Y$
over $K$, and give $L=K((s))$ the $s$-adic valuation and the
coefficientwise derivations. The point $a=(s,s^2)$ has residue
$b=(0,0)$. At $a$ the translated pair is
$((s+T_1)^2,s^2+T_2)$, whose degree-zero pair $(s^2,s^2)$
represents $[1:1]$. At $b$ the translated pair is $(T_1^2,T_2)$,
whose lowest homogeneous pair is $(0,T_2)$ in degree one. Hence
\[
       \res_v\bigl(\rho(a;T)\bigr)=[1:1],
       \qquad \rho(b;T)=[0:1].
\]

In this example every graph equation specializes identically to zero
at $b$. Indeed, let $z$ be a differential indeterminate over $K$ and
work in $K\langle z\rangle((s))$, with $s$ constant under all
derivations. Substituting $(X,Y)=(sz,s^2z)$ gives $X^2/Y=z$, so
$P(sz,s^2z,z)=0$ for every $P\in J$. Taking constant terms in $s$
gives $P(0,0,z)=0$. Since $z$ is differentially transcendental
over $K$, this says that $P(0,0,Z)=0$ in $K\{Z\}$.
\end{remark}

\section{The Proof of the Conjecture}\label{sec:uniqueness}

 Write $A=K\{Y,Z\}/J$.
The pair $(b,z)$ in a differential field extension of $K$  is called a \emph{graph point} if $P(b,z)=0$ for every $P\in J$.
If each $P(b,Z)$ is the zero differential polynomial, every $z$
gives a graph point over $b$. If some $P(b,Z)$ is nonzero, the
following proposition determines the ordinate of every such point in this chart.

\begin{proposition}\label{prop:uniqueness}
Let $(b,z)$ be a graph point in a differential field extension of
$K$. If $P(b,Z)\ne0$ for some $P\in J$, then $\rho(b;T)=[z:1]$.
\end{proposition}

\begin{proof}
  Let $E$ be a differential field extension of $K$ containing $(b,z)$.
Since $(b,z)$ is a graph point, evaluation induces a differential
$K$-algebra homomorphism
\[
 e_{b,z}:A\longrightarrow E,\qquad H+J\longmapsto H(b,z).
\]
Its kernel $\pp_0$ is a prime differential ideal. Fix $P\in J$
with $P(b,Z)\ne0$ and choose a coefficient $C\in K\{Y\}$ of a differential
monomial in $Z$ such that $C(b)\ne0$. For each prime
differential ideal $\pp\subseteq\pp_0$, write
$\kappa(\pp)=\Frac(A/\pp)$ and let $(a_\pp,z_\pp)$ be the images of
$(Y,Z)$ in $\kappa(\pp)$. Since the class of $C$ lies outside
$\pp_0$, it also lies outside $\pp$. Thus $C(a_\pp)\ne0$, and
$P(a_\pp,Z)$ is a nonzero differential polynomial over $\kappa(\pp)$.

Let $\mathcal B$ be the set of prime differential ideals
$\pp\subseteq\pp_0$ for which $\rho(a_\pp;T)\ne[z_\pp:1]$.
Since $\pp_0=\ker e_{b,z}$, evaluation identifies $A/\pp_0$ with
the differential $K$-algebra generated by $b,z$ in $E$. On passing
to fraction fields, this identifies $\kappa(\pp_0)$ with
$K\langle b,z\rangle$ and sends $(a_{\pp_0},z_{\pp_0})$ to $(b,z)$.
Under this identification the expansions in~\eqref{eq:expansion}
at $a_{\pp_0}$ and at $b$ have the same coefficients, hence the
same first nonzero homogeneous pair. Consequently
\[
       \rho(b;T)\ne[z:1]\quad\Longleftrightarrow\quad
       \pp_0\in\mathcal B.
\]
It therefore suffices to show that $\mathcal B$ is empty. We consider
all primes contained in $\pp_0$ so that, if a counterexample exists,
we can choose a minimal one and apply the valuative comparison from
a strictly smaller prime.

 Write $\mathcal M$ for the monomials in~\eqref{eq:expansion} and
identify $f_M,g_M,Z$ with their images in $A$. At $\pp$ the least
degree is
\[
 d_\pp=\min\{d_M:f_M\notin\pp\ \text{or}\ g_M\notin\pp\}.
\]
Comparing coefficients in degree $d_\pp$ gives
\begin{equation}\label{eq:mismatch}
 \pp\in\mathcal B\quad\Longleftrightarrow\quad
 f_M-Zg_M\notin\pp\ \text{for some }M\text{ with }d_M=d_\pp .
\end{equation}
Thus membership in $\mathcal B$ is determined by the intersection
of $\pp$ with the finite set
\[
 S=\bigcup_{M\in\mathcal M}\{f_M,g_M,f_M-Zg_M\}.
\]
Indeed, $f_M$ and $g_M$ determine which coefficient pairs survive,
and $f_M-Zg_M$ tests proportionality to $(z_\pp,1)$ in the first
degree in which a pair survives. 

  Suppose $\mathcal B$ is nonempty. For any nonempty chain
$\mathcal C\subseteq\mathcal B$, choose $\pp_*\in\mathcal C$
minimizing $|\pp_*\cap S|$. The sets $\pp\cap S$, for
$\pp\in\mathcal C$, are nested and finite, so $\pp_*\cap S$ is
contained in all of them. The intersection
$\mathfrak r=\bigcap_{\pp\in\mathcal C}\pp$ is again prime and
differential, and
\[
 \mathfrak r\cap S=\bigcap_{\pp\in\mathcal C}(\pp\cap S)=\pp_*\cap S.
\]
Thus $\mathfrak r\in\mathcal B$ is a lower bound for $\mathcal C$,
and Zorn's lemma gives an inclusion-minimal member $\pp$ of
$\mathcal B$. The zero prime lies outside $\mathcal B$, since the
degree-zero coefficient pair at the generic point is
$(F(t),G(t))=G(t)(u,1)$ with $G(t)\ne0$. Hence $\pp\ne(0)$. 

 By Theorem~\ref{lem:valuation} there is a prime differential ideal
$\qq\subsetneq\pp$ and a monotone real valuation $v$ on $\kappa(\qq)$
whose valuation ring contains $A/\qq$ and whose maximal ideal
contracts to $\pp/\qq$. Reduction on $A/\qq$ induces a differential
embedding $\kappa(\pp)\hookrightarrow k_v$ under which
$(a_\pp,z_\pp)$ is the residue of $(a_\qq,z_\qq)$. In particular,
$v$ is trivial on $K$. By minimality, $\qq\notin\mathcal B$, so
$\rho(a_\qq;T)=[z_\qq:1]$. The nonzero coefficient $C(a_\pp)$
remains nonzero in $k_v$. Proposition~\ref{lem:comparison}
therefore gives
\[
 \rho(a_\pp;T)=\res_v[z_\qq:1]=[z_\pp:1],
\]
where the last equality uses $z_\qq\in\OO_v$. This equality holds
already over $\kappa(\pp)$ by injectivity, contradicting
$\pp\in\mathcal B$. Hence $\mathcal B$ is empty. 
\end{proof}

 We now pass from uniqueness in one chart to the projective fiber.
The two graph ideals give the same test for fullness. When the fiber
is proper, the constant-parameter substitution produces an admitted
point over $K\langle b\rangle$, and Proposition~\ref{prop:uniqueness}
in the two charts excludes every other point. 

\begin{proof}[{\bf Proof of Theorem~\ref{thm:main}}]
 Fix $b\in\U^n$ and put $E=K\langle b\rangle$. Let $T$ be a tuple
of differential indeterminates over $E$. Write $d=d(b)$, and let
$F_d(T)$ and $G_d(T)$ be the homogeneous parts of degree $d$ of
$F(b+T)$ and $G(b+T)$, respectively. Thus
$\rho(b;T)=[F_d(T):G_d(T)]$.  
 We must show that either the
fiber over $b$ is all of $\PP^1(\U)$, or it consists of the single
point $\rho(b;T)$ and this point belongs to $\PP^1(E)$.

 We claim that the   graph ideals $J\subset K\{Y,Z\}$ and
$J^\vee\subset K\{Y,W\}$ satisfy
\begin{equation}\label{eq:relation}
 \bigl(\exists P\in J,\ P(b,Z)\ne0\bigr)
 \quad\Longleftrightarrow\quad
 \bigl(\exists Q\in J^\vee,\ Q(b,W)\ne0\bigr).
\end{equation}
Substitution $Z=1/W$ induces a differential isomorphism
$K\{Y,Z\}[1/Z]\simeq K\{Y,W\}[1/W]$ with inverse $W=1/Z$.
Given $P\in J$ with $P(b,Z)\ne0$, choose $e\ge0$ such that
$Q=W^eP(Y,1/W)$ belongs to $K\{Y,W\}$. Evaluation at $(t,1/u)$
gives $Q(t,1/u)=0$, so $Q\in J^\vee$. After specializing $Y$ to
$b$, inversion still induces an isomorphism
$E\{Z\}[1/Z]\simeq E\{W\}[1/W]$. Hence $P(b,1/W)\ne0$ and
$Q(b,W)\ne0$. Interchanging $Z$ and $W$ gives the converse. 

\vskip3pt
 If neither side of~\eqref{eq:relation} holds, every graph equation
in either chart specializes to the zero differential polynomial.
By~\eqref{eq:fiber-charts}, every point of $\PP^1(\U)$ is then
admitted. 

 Suppose now that the equivalent conditions in~\eqref{eq:relation}
hold. At least one of $F_d(T)$ and $G_d(T)$ is nonzero. Interchange
the charts, if necessary, so that $G_d(T)\ne0$. Put
$L=E\langle T\rangle$ and adjoin an indeterminate $s$ with
$\delta_i s=0$ for every $i$. Substituting $sT$ for $T$ multiplies
a differential monomial of total degree $j$ by $s^j$. Hence
\[
 F(b+sT)=s^d U(s),\qquad G(b+sT)=s^d V(s),
\]
where $U,V\in L[s]$, $U(0)=F_d(T)$ and $V(0)=G_d(T)\ne0$. Thus
the quotient $u_s=U(s)/V(s)$ is regular at $s=0$, that is,
$u_s\in R=L[s]_{(s)}$, since its
presentation $U(s)/V(s)$ has denominator nonzero at zero. Evaluation
on the local ring $R=L[s]_{(s)}$ of such fractions is given by
\[
 \varepsilon:R\longrightarrow L,\qquad
 p(s)/q(s)\longmapsto p(0)/q(0)
 \qquad(q(0)\ne0)
\]
which commutes with each derivation. Indeed, $\delta_i$ acts
coefficientwise on $L[s]$, and the quotient rule introduces only a
power of $q$ in the denominator. 

 Let $r=\varepsilon(u_s)=F_d(T)/G_d(T)$. For any $H\in J$, the
graph identity $H(Y,F/G)=0$ gives $H(b+sT,u_s)=0$ since
$G(b+sT)\ne0$. All coordinates of $(b+sT,u_s)$ lie in $R$, so
applying $\varepsilon$ yields $H(b,r)=0$. In particular, for
$P\in J$ as in~\eqref{eq:relation}, the nonzero differential
polynomial $P(b,Z)\in E\{Z\}$ vanishes at $r$. Thus $r$ is
differentially algebraic over $E$. Since $r\in E\langle T\rangle$,
Lemma~\ref{lem:free} gives $r\in E$. Consequently $(b,r)$ is a
graph point in $\U$, and $r$ is admitted in the chosen chart. 

  It remains to prove uniqueness. For each finite admitted value $z$ in this chart,
Proposition~\ref{prop:uniqueness} gives $[z:1]=[r:1]$, that is, $z=r$. If the
point at infinity in this chart were admitted, $(b,0)$ would be a
graph point in the reciprocal chart. By~\eqref{eq:relation} the same
proposition would give $[G_d(T):F_d(T)]=[0:1]$, contradicting
$G_d(T)\ne0$. This proves that the fiber is a singleton. In the
original homogeneous coordinates its unique point is
$[F_d(T):G_d(T)]$, which belongs to $\mathbb P^1(E)$. 
\end{proof}

\begin{proof}[Proof of Corollary~\ref{cor:ritt}]
 Take $b=0$ in Theorem~\ref{thm:main}. Then
$K\langle b\rangle=K$, and \eqref{eq:fiber-charts} identifies
the points of the fiber with the values admitted at the origin,
including infinity. 
\end{proof}

\begin{example} 
 Let $X,Y$ be ordinary differential indeterminates over
$K$ and put $u=XY'/Y$. Under the substitution $(X,Y)\mapsto(sX,sY)$,
with $s$ a constant indeterminate, the quotient becomes $sXY'/Y$.
Specializing $s$ to zero gives the admitted value $0$ at the origin.
To obtain a second value, first keep $X$ fixed and work over
$E=K\langle X\rangle$. As polynomials in $Y$ and its derivatives,
$XY'$ and $Y$ are homogeneous of degree one, and their ratio lies
outside $E$. By Blum's observation following Proposition~1
\cite[\S1, p.~678]{blum}, the specialization $Y\mapsto0$ therefore
admits the value $1$ over $E$. Since $X$ is differentially
transcendental over $K$, we may then specialize $X$ to zero over $K$.
Thus
\[
 (X,Y,u)\rightsquigarrow_K(X,0,1)\rightsquigarrow_K(0,0,1).
\]
The fiber at the origin contains both $0$ and $1$, so
Theorem~\ref{thm:main} gives the whole projective line. 

 A monotone valuation on $K\langle X,Y\rangle$ specializing
$(X,Y,u)$ to $(0,0,1)$ would give $v(X)>0$ and $v(u)=0$.
Monotonicity, however, gives $v(Y')\ge v(Y)$, whence
\[
 v(u)=v(X)+v(Y')-v(Y)\ge v(X)>0 .
\]
Let $\pp$ be the kernel of evaluation at $(0,0,1)$ in the graph
coordinate ring $A$. Taking $\qq=(0)$ in Theorem~\ref{lem:valuation}
would give a monotone valuation on $\Frac(A)=K\langle X,Y\rangle$
with center $\pp$, realizing the specialization just ruled out.
Thus for this $\pp$ one must have $(0)\subsetneq\qq\subsetneq\pp$. 
\end{example}

 We conclude with the example concerning Blum's order hypothesis
mentioned in the introduction. 

\begin{remark} 
\label{rem:blum-order}
  In the ordinary differential case, let $F,G\in K\{Y\}$ be
nonzero and relatively prime, with $F(0)=G(0)=0$.
Blum proved that the fiber at zero is full if $[Y]$ is
not a component of $ZG-F$ and $\operatorname{ord}(FG)=1$,
and asked whether the order hypothesis could be omitted
\cite[Theorem~1(iii), p.~679]{blum}.
Here a component is a prime differential ideal minimal over
the differential ideal generated by the equation. 

 Substituting $Y'$ for $Y$ in Ritt's example $(Y')^2/Y$
\cite[p.~681]{blum} gives the quotient and graph equation
\[
 f(Y)=\frac{(Y'')^2}{Y'},\qquad H(Y,Z)=ZY'-(Y'')^2.
\]
The numerator and denominator are relatively prime and their
product has order two. Since $K\{Y,Z\}/[Y']\simeq K\{Z\}[Y]$
is a domain, $[Y']$ is prime. The inclusions
$H\in[Y']\subsetneq[Y]$ therefore show that $[Y]$ is not a
component of $H$. 

To compute the fiber, put $x=Y'$, $z=(x')^2/x$ and $w=1/z$.
Differentiation gives the graph identities
\[
 z^2=2x''z-x'z',\qquad 1=2x''w+x'w'.
\]
Under the differential specialization $Y\mapsto0$, both $x'$ and
$x''$ vanish. The first identity forces every finite ordinate to be
zero, and the second excludes infinity. Projective nonemptiness
\cite[Prop.~1, p.~678]{blum} now gives $\operatorname{Fib}_0(f)=\{0\}$.
Thus the order restriction in Blum's criterion is necessary. 
\end{remark}


\begin{thebibliography}{10}
\bibitem{blum}
P. Blum, \emph{Rational functions on differentially closed sets},
Amer. J. Math. \textbf{94} (1972), 676--684.
\href{https://doi.org/10.2307/2373750}{doi:10.2307/2373750}.
\bibitem{engler-prestel}
A. J. Engler and A. Prestel, \emph{Valued Fields},
Springer Monographs in Mathematics, Springer, Berlin, 2005.
\bibitem{ffh}
J. Fang, J. Freitag and M. Harrison-Trainor,
\emph{The Ritt problem is decidable}, preprint (2026),
\href{https://arxiv.org/abs/2611.27008}{arXiv:2611.27008}.
\bibitem{icm}
E. R. Kolchin, \emph{Some problems in differential algebra},
Proceedings of the International Congress of Mathematicians
(Moscow, 1966), Mir, Moscow, 1968, 269--276.
\bibitem{kolchin}
E. R. Kolchin, \emph{Differential Algebra and Algebraic Groups},
Pure and Applied Mathematics, vol.~54, Academic Press, New York, 1973.

\bibitem{morrison}
 S. D. Morrison, \emph{Differential specializations in integral
and algebraic extensions}, Amer. J. Math. \textbf{101} (1979),
1381--1399. \href{https://doi.org/10.2307/2374147}{doi:10.2307/2374147}. 
\bibitem{ng}
G. Ng, \emph{Taylor morphisms}, J. Algebra \textbf{701} (2026),
219--255.
\href{https://doi.org/10.1016/j.jalgebra.2026.04.028}{doi:10.1016/\allowbreak j.jalgebra.\allowbreak 2026.04.028}.
\bibitem{ritt1932}
J. F. Ritt, \emph{Differential Equations from the Algebraic Standpoint},
American Mathematical Society Colloquium Publications, vol.~14,
American Mathematical Society, New York, 1932.

\bibitem{ritt}
J. F. Ritt, \emph{Indeterminate expressions involving a function
and its derivatives}, Monatsh. Math. Phys. \textbf{43} (1936), 97--104.
\href{https://doi.org/10.1007/BF01707591}{doi:10.1007/BF01707591}.

\bibitem{ritt1950}
J.~F. Ritt,
\emph{Differential Algebra},
American Mathematical Society Colloquium Publications,
vol.~33, American Mathematical Society,
New York, 1950.

\end{thebibliography}
\end{document}